\documentclass[12pt]{article}
\usepackage{amsmath,amssymb}
\usepackage{amsthm}
\usepackage{amsfonts}   
\author{Simona Settepanella}
\title{\textbf{A stability-like theorem for cohomology of Pure Braid Groups of 
the series $A$, $B$ and $D$}}

\DeclareMathOperator{\cl}{cl}

\newcommand{\Q}{\mathbb Q}
\newcommand{\N}{\mathbb N}
\newcommand{\Z}{\mathbb Z}
\newcommand{\R}{\mathbb R}
\newcommand{\C}{\mathbb C}
\newtheorem{thm}{Theorem}[section]
\newtheorem{prop}{Proposition}[section]

\newtheorem{lemma}{Lemma}[section]
\newtheorem{rem}{Remark}[section]

\newcommand{\bW}{\mathbf W}
\newcommand{\bQ}{\mathbf Q}
\newcommand{\bY}{\mathbf Y}
\newcommand{\bX}{\mathbf X}
\newcommand{\aenne}{\mathbf{A_n}}
\newcommand{\benne}{\mathbf{B_n}}
\newcommand{\denne}{\mathbf{D_n}}

\newcommand{\lgr}{\longrightarrow}
\newcommand{\ver}{\mathit  v}
\newcommand{\s  }{\mathit  s}

\newcommand{\A}{\mathcal A}

\newcommand{\st}{\mathcal S}
\newcommand{\As}{{\mathbb A}^n}
\newcommand{\aemme}{\mathbf{A_n}}
\newcommand{\bemme}{\mathbf{B_n}}
\newcommand{\demme}{\mathbf{D_n}}

\begin{document}

\maketitle

{\small ABSTRACT. Consider the ring $R:=\Q[\tau,\tau^{-1}]$ of Laurent 
polynomials in the variable $\tau$. The Artin's Pure Braid Groups (or 
Generalized 
Pure Braid Groups) act over $R,$ where the action of every
standard generator is the multiplication by $\tau$. 
In this paper we consider the cohomology of such groups with coefficients in 
the module $R$ (it is well known that such cohomology is strictly related to 
the untwisted integral cohomology of the Milnor fibration naturally associated
to the reflection arrangement). We give a sort of \textit{stability}
theorem for the cohomologies of the infinite series $A$, $B$ and $D,$   finding
that these cohomologies stabilize, with respect to the natural inclusion, 
at some number of copies of the trivial $R$-module $\Q$. We also give a 
formula which computes this number of copies. }

\section{Introduction}

Let $(\bW,S)$ be a finite Coxeter system realized as a reflection group 
in $\R^n$,
$\mathcal{A}(\bW)$ the 
arrangement in $\C^n$ obtained by complexifying the reflection hyperplanes
of $\bW$.Let
$$ 
\bY(\bW)=\bY(\mathcal{A}(\bW))=\C^n \setminus \cup_{H \in \mathcal{A}(\bW)}H.
$$ 
be the complement to the arrangement,
then $\bW$ acts freely on $\bY(\bW)$ and the fundamental group $G_W$ of 
the orbit space $\bY(\bW)/\bW$ is the so called \emph{Artin group} associated
to $\bW$ (see \cite{Bou}). Likewise the fundamental group $P_W$ of $\bY(\bW)$ 
is the \emph{Pure Artin group} or the pure braid group of the series $\bW$.
It is well known (\cite{Brie}) that these spaces $\bY(\bW)$ $(\bY(\bW)/\bW)$
are of type $K(\pi,1)$, so there cohomologies equal that of $P_W$ $(G_W)$.

The integer cohomology of $\bY(\bW)$ is well known (see \cite{Brie},\cite{OS},
\cite{b1},\cite{b5})
and so is the integer cohomology of the Artin groups associated to finite
Coxeter groups (see \cite{b9},\cite{b8},\cite{boss1}).

Let $R=\Q[\tau,\tau^{-1}]$ be the ring of rational Laurent polynomials.
The R can be given a structure of module over the Artin group $G_W$, where
standard generators of $G_W$ act as $\tau$-multiplication.
  
In \cite{b10} and \cite{b13} the authors compute the cohomology of all Artin 
groups associated to finite Coxeter groups with coefficients in the previous
module.

In a similar way we define a $P_W$-module $R_{\tau}$, where standard generators
of $P_W$ act over the ring $R$ as $\tau$-multiplication.

Equivalently, one defines an abelian local system (also called $R_{\tau}$) 
over $\bY(\bW)$ with fiber $R$ and local monodromy around each hyperplane
given by $\tau$-multiplication 
(for local systems on $\bY(\bW)$ see \cite{b14},\cite{b15}).
   
In this paper we are going to consider the cohomology of $\bY(\bW)$ with 
local coefficients $R_{\tau}$, for the finite Coxeter groups of the series  
$A$, $B$ and $D$ (see \cite{Bou})
(that is equivalent to the cohomology of $P_W$ 
with coefficients in $R_{\tau}$).

Our aim is to give a sort of \emph{``stability''} theorem for these 
cohomologies
(for stability in the case of Artin groups see \cite{Decsal}) .

Denote by $\varphi_i$ the i-th cyclotomic polynomial and let be 
$$
\{\varphi_i\}:=\Q[\tau,\tau^{-1}]/(\varphi_i)=\Q[\tau]/(\varphi_i)
$$
thought as $R$-module.By its definition $\{\varphi_1\}=1-\tau$ so that
$\{\varphi_1\}=\Q$. 

Notice that by identification $\Q[\tau,\tau^{-1}] \cong \Q[\Z]$, the sums 
of copies of $\{\varphi_1\}$ are the unique trivial $\Z$-modules. We obtain

\begin{thm}\label{teor:stab} Let  
$\bW$ be a Coxeter group of type $A_n$, then for $n \geq 3k-2$ the cohomology 
group $H^k(\bY(A_n),R_{\tau})$
is a trivial $\Z$-module.\\
Analog statement holds for $\bW$ of type $B_n$ in the rang $n \geq 2k-1$ 
and for $\bW$ of type $D_n$ in the rang $n \geq 3k-1$.
\end{thm}

The proof of this theorem is obtained extending the 
methods developed in \cite{b10} and using some known results about the
global Milnor fibre $F(\bW)$ of the complement $\bY(\bW)$.

We recall briefly that if $H \in \mathcal{A}=\mathcal{A}(\bW) $ and 
$\alpha_H \in \C[x_1,\ldots ,x_n]$ is a linear form s.t. $H=ker(\alpha_H)$, 
then the global
Milnor fibre $F(\bW)$ is a complex manifold of dimension $n-1$ given by 
$F(\bW)=Q^{-1}(1)$ where $Q=Q(\mathcal{A})=\prod_{H \in \mathcal{A}}\alpha_H$
is the \emph{defining polynomial} for $\mathcal{A}$.

It is well known (see also \cite{Denham}) that, over $R$, there is a 
decomposition  
\begin{equation*}
H^*(F(\bW),\Q) \simeq \bigoplus_{i \mid \sharp (\mathcal{A}(\bW))}
(R/(\varphi_i))^{\alpha_i}=
\bigoplus_{i \mid \sharp (\mathcal{A}(\bW))}\{\varphi_i\}^{\alpha_i}.
\end{equation*}
the action on the left being that induced by monodromy.

Since $F(\bW)$ is homotopy-equivalent to an infinite cyclic cover of 
$\bY(\bW)$, there is an isomorphism of $R$-modules 
$$
H^*(F(\bW),\Q) \simeq H^*(\bY(\bW),R_{\tau})
$$
and then 
\begin{equation}\label{spezz}
H^*(\bY(\bW),R_{\tau}) \simeq \bigoplus_{i \mid \sharp (\mathcal{A}(\bW))}
\{\varphi_i\}^{\alpha_i}.
\end{equation}
The other tool we use is a suitable filtration by subcomplexes of the 
algebraic Salvetti's CW-complex ($C(\bW)$,$\delta$) coming from \cite{boss0} 
(see also \cite{boss2}, \cite{boss1}), which we recall in the next paragraph. 

Finally we use the universal coefficients theorem to compute the dimensions
of the above cohomologies as vector spaces over the rationals.

\begin{thm}\label{teor:esp} 
In the range specified in theorem \ref{teor:stab} one has:
$$
\mathtt{rk }H^{k+1}(Y(\bW),R_{\tau})=\sum_{i=0}^{k}(-1)^{(k-i)}
\mathtt{rk }H^i(Y(\bW),\Z).
$$ 
\end{thm}
So one reduces to compute the dimensions of
the Orlik-Solomon algebras of $\A(\aenne)$, $\A(\benne)$ and $\A(\denne)$ 
(see \cite{b11}).

\section{Salvetti's Complex}

Let $\bW$ be a finite group generated by reflections in the affine space $\As(\R )$.
Let $\overline{\A}(\bW)=\{H_j\}_{j \in J}$ be the arrangement in $\As$ defined by 
the reflection 
hyperplanes of $\bW$. We need to recall briefly some notations and results 
from \cite{boss0} for the particular case of Coxeter arrangements. 
$\overline{\A}(\bW)$ induces a stratification
$\st=\st(\bW)$ of $\As$
into facets (see \cite{Bou}). The set $\st$ is partially
ordered by $F > F^{\prime}$ iff $ F' \subset \cl (F)$. 
We shall indicate by $\bQ=\bQ(\bW)$ the cellular complex which is
\textit{dual} to $\st$. 
In a standard way, this can be realized inside $\As$ by baricentrical 
subdivision of the facets: inside each codimension j facet
$F^j$ of $\st$ choose one point $\ver( F^j )$ and consider
the simplexes
$$
\s(F^{i_{0}},\cdots,F^{i_{j}}) = \{ \sum_{k=0}^{j}\lambda_k \ver(F^{i_{k}}) :
\sum_{k=0}^{j}\lambda_k=1, \lambda_k \in [0,1]\}
$$
where $F^{i_{k+1}} < F^{i_{k}}, k=0,\cdots,j-1$. The dimension j cell 
$e^j(\overline{F}^j)$ which is dual to $\overline{F}^j$ is obtained by
taking  the union
$$
\cup \s(F^0,\cdots ,F^{j-1},\overline{F}^j)
$$
over all chains $\overline{F}^j < F^{j_1} < \cdots < F^0$.
So $\bQ=\cup e^j(F^j)$, the union being over all facets of $\st$.

One can think of the $1-skeleton$ $\bQ_1$ as a graph (with vertex-set 
the $0-skeleton$ $ \bQ_0$) and can define the combinatorial distance
between two vertices $\ver,\ver^{\prime}$ as the minimum number of edges in an
edge-path connecting $\ver$ and $\ver^{\prime}$.

For each cell $e^j$ of $\bQ$ one indicates by $V(e^j)=\bQ_0 \cap e^j$ the
$0$-\textit{skeleton} of $e^j$. One has
\begin{prop}\label{length2}
Given a vertex v $\in \bQ_0$ and a cell $e^i \in \bQ$, there is
a unique vertex $\underline{w}(v,e^i) \in V(e^i)$ with the lowest
combinatorial distance from v, i.e.:
\begin{equation*}
d(v,\underline{w}(v,e^i)) < d(v,v^{\prime}) \textrm{ if } v' \in V(e^i) \setminus
\{\underline{w}(v,e^i)\}.
\end{equation*}
If $e^j \subset e^i$ then $\underline{w}(v,e^j) =
\underline{w}(\underline{w}(v,e^i),e^j)$.
\end{prop} 

Let now $\A(\bW)$ denote the \textit{complexification} of $\overline{\A}(\bW)$,
and 
$\bY(\bW)=\C^n \setminus \cup_{j\in J}H_{j,\C}$ the complement of the
complexified arrangement. Then $\bY(\bW)$  is homotopy
equivalent to the complex $\bX(\bW)$ which is constructed as follows 
(see \cite{boss0}).

Take a cell $e^j=e^j(F^j)=\cup \s(F^0,\cdots,F^{j-1},F^j)$ of $\bQ$ as
defined above and let $v \in V(e^j)$. Embed each simplex
$\s(F^0,\cdots,F^j)$ into $\C^n$ by the formula 
\begin{equation}\label{immersione}
\begin{split}
\phi_{v,e_j} (\sum_{k=0}^{j}\lambda_k & \ver (F^k))=\\
&\sum_{k=0}^{j}\lambda_k \ver (F^k) + i \sum_{k=0}^{j}\lambda_k
(\underline{w}(v,e^k)-\ver (F^k)).
\end{split}
\end{equation}
It is shown in \cite{boss0} (see also \cite{boss1}):

\textit{(i)} the preceding formula defines an embedding of $e^j$ into 
$\bY(\bW)$;

\textit{(ii)} if $E^j = E^j(v,e^j)$ is its image, then varying $e^j$ and $\ver$
one obtains 

a cellular complex
$$
\bX(\bW) = \cup E^j
$$

which is homotopy equivalent to $\bY(\bW)$.\\
The previous result allows us to make cohomological computations over 
$\bY(\bW)$ by using the complex $\bX(\bW)$.

In \cite{boss1} (see also \cite{boss3}) 
the authors give a new combinatorial description of the
stratification $\st$ where the action of $\bW$ is more
explicit. They prove that if 
$S$ is the set of reflections with respect to the walls of the fixed base
chamber $C_0$, then a cell in $\bX(\bW)$ is of the form $E=E(w,\Gamma)$ with 
$\Gamma \subset S$ and $w \in 
\bW$. The action of $\bW$ is written as 
\begin{equation}\label{prodperaz}
\sigma.E(w,\Gamma)=E(\sigma w,\Gamma),
\end{equation}
where the factor $\sigma .w$ is just multiplication in $\bW$.\\

We prefer at the moment to deal with chain complexes and boundary operator 
coming from $\bX(\bW)$ instead of cochain and coboundary. Then we will deduce
cohomological results by standard methods.  

We define a rank-1 local system on $\bY(\bW)$ with coefficients in an 
unitary ring $A$  by assigning an unit
$\tau_j=\tau(H_j)$ (thought as a multiplicative operator) to each hyperplane
$H_j \in \A$.
Call $\overline{\tau}$ the collection of $\tau_j$ and 
$\mathcal{L}_{\overline{\tau}}$ the 
corresponding local system.
Let $C(\bW,\mathcal{L}_{\overline{\tau}})$ be the free graduated $A$-module 
with basis all $E(w,\Gamma)$.

We use the natural identification between the elements of the group and the 
vertices of $\bQ_0$, given by $w \leftrightarrow w.v_0$. Here $v_o \in \bQ_0$
is contained in the fixed base chamber $C_0$.

Then $u(w, w^{\prime})$ will denote the \emph{``minimal positive path''} 
joining the 
corresponding vertices $v$ and $v^{\prime}$ in the $1$-skeleton $\bX(\bW)_1$
of $\bX(\bW)$ (see \cite{boss0}).

The local system $\mathcal{L}_{\overline{\tau}}$ defines for each edge-path 
$c$ in $\bX(\bW)_1$, $c:w \rightarrow w^{\prime}$ an isomorphism 
$c_*:A \rightarrow A$ such that for all $d:w \rightarrow w^{\prime}$ homotopic 
to $c$, $c_*=d_*$ and 
for all $f:w^{\prime \prime} \rightarrow w$, $(cf)_*=c_*f_*$.

Then the set $\{s_0(w).E(w,\Gamma)\}_{|\Gamma|=k}$, where
$s_0(w):=u(1,w)_*(1)$, is a linear basis of 
$C_{k}(\bW,\mathcal{L}_{\overline{\tau}})$. 

Let now $T=\{wsw^{-1} | s \in S, w \in \bW\}$, the set of reflections in $\bW$
and 
$$
\overline{\bW }=\{\mathbf{s}(w)=(s_{i_1},\cdots ,s_{i_q}) | 
w = s_{i_1}\cdots s_{i_q} \in \bW\},
$$
then for each $\mathbf{s}(w) \in \overline{\bW }$ and $t \in T$, we set

i) $\Psi(\mathbf{s}(w))=(t_{i_1},\cdots,t_{i_q})$ with $t_{i_j}=
(s_{i_1}\cdots s_{i_{j-1}})s_{i_j}(s_{i_1}\cdots s_{i_{j-1}})^{-1} \in T$ 

ii) $\overline{\Psi(\mathbf{s}(w))}=\{t_{i_1},\cdots,t_{i_q}\}$ 
 
iii)$\eta(w,t)=(-1)^{n(\mathbf{s}(w),t)}$ with 
$n(\mathbf{s}(w),t)=\sharp\{j | 1 \leq j \leq q \textnormal{ and } 
t_{i_j}=t\}$.\\     
Moreover if $t \in T$ is the reflection relative to the hyperplane $H$,
then we set $\tau(t)=\tau(H)$.

We define
\begin{equation}\label{bordo}
\begin{split}
\partial_k(\s_0(w).&E(w,\Gamma)) =\\ 
&\sum_{\sigma\in \Gamma}
\sum_{\beta\in\bW^{\Gamma\setminus\{\sigma\}}_{\Gamma}}(-1)^{l(\beta)+\mu
(\Gamma,\sigma)}\tau(w,\beta)\s_0(w\beta).
E(w\beta,\Gamma\setminus\{ \sigma \}).
\end{split}
\end{equation}
where $\tau(w,\beta )=\prod\limits_{\substack{\scriptstyle t \in \overline{\Psi(\mathbf{s}(w))} \\ \scriptstyle \eta(w,t)=1}} \tau(t)$,
 and $\mu(\Gamma,\sigma)=\sharp\{i \in \Gamma | i \leq \sigma\}$.

We have the following (see \cite{boss3}, \cite{Simona})

\begin{thm}
$H_*(C(\bW),\mathcal{L}_{\overline{\tau}})\cong 
H_*(C(\bW,\mathcal{L}_{\overline{\tau}}),\partial)$.
\end{thm}

We have a similar result for the cohomology.

\section{A filtration for the complex $(C(\bW),\partial)$}

Let $(\bW,S)$ be a finite Coxeter system with $S=\{s_1,\cdots ,s_n\}$.
We are interested in the cohomology of $C(\bW)$ 
(equivalently $\bY(\bW)$) with coefficients in $R_{\tau}$ (see introduction).

In this case the boundary operator defined in (\ref{bordo}) becomes
\begin{equation}\label{rel:bordoR}
\partial(E(w,\Gamma)) = 
\sum_{\sigma\in \Gamma}
\sum_{\beta\in\bW^{\Gamma\setminus\{\sigma\}}_{\Gamma}}(-1)^{l(\beta)+\mu
(\Gamma,\sigma)}\tau^{\frac{l(\beta)+l(w)-l(w\beta)}{2}}
E(w\beta,\Gamma\setminus\{ \sigma \})
\end{equation}
where $\tau$ is the variable in the ring $R$. 

From (\ref{spezz}) and universal coefficients theorem it follows that 
\begin{equation}\label{comom}
H^*(C(\bW),R_{\tau})=H_{*-1}(C(\bW),R_{\tau}).
\end{equation}

For each integer $0 \leq k \leq n$ denote by 
$S_k=\{s_1,\cdots ,s_k\} \subset S$ and $S^k=S \setminus S_k$.
We define the graduated $R$-submodules of $C(\bW)$:
\begin{equation*}
\begin{split}
&G_n^k(\bW):=\mathop{\sum_{w \in \bW}}_{\Gamma \subset S_k} R.E(w,\Gamma) \\
&F_n^k(\bW):=\mathop{\sum_{w \in \bW}}_{\Gamma \supset S^{n-k}} R.E(w,\Gamma).
\end{split}
\end{equation*}

There is an obvious inclusion
\begin{equation} \label{eqn3}
i_{n,h} : G_n^{n-h}(\bW) \lgr G_n^n(\bW)=C(\bW).
\end{equation}

Each $G_n^k(\bW)$ is preserved by the induced boundary map and we get a 
filtration by subcomplexes of $C(\bW)$:
$$
C(\bW)=G_n^n(\bW) \supset G_n^{n-1}(\bW) \cdots \supset  G_n^1(\bW) \supset G_n^0(\bW).
$$
The quotient module $G_n^n(\bW)/G_n^{n-1}(\bW)$ is exactly $F^1_n(\bW)$ 
which becomes an algebraic complex with the induced 
boundary map.

We give iteratively to $F^k_n(\bW)$, $k \geq 2$, a structure of complex by 
identifying it with the cokernel of the map:
\begin{equation*}
\begin{split}
&i_n[k] : G_n^{n-(k+1)}(\bW)[k] 
\lgr F^k_n(\bW),\\
&i(E(w,\Gamma)))=E(w,\Gamma \cup S^{n-k}).
\end{split} 
\end{equation*}

Here $M[k]$ denotes, as usual, $k$-augmentation of a complex $M$; so $i_n[k]$
is degree preserving.

By construction $i_n[k]$ gives rise to the exact sequence of complexes
\begin{equation} \label{eqn:suc.es.}
0 \lgr G_n^{n-(k+1)}(\bW)[k] \lgr 
F^k_n(\bW) \lgr F^{k+1}_{n}(\bW) \lgr 0.
\end{equation}

Let $\Gamma \subset S$ and let $\bW_{\Gamma}$ be the \emph{parabolic subgroup} 
of $\bW$ generated by $\Gamma$. Recall from \cite{Bou} the following
\begin{prop}
Let $(\bW,S)$ be a Coxeter system. Let $\Gamma\subset S$. The following statements hold.

(i) $(\bW_\Gamma, \Gamma )$ is a Coxeter system.

(ii) Viewing  $\bW_\Gamma$ as a Coxeter group with length function $\ell_\Gamma$, $\ell_S=\ell_\Gamma$ on $\bW_\Gamma$.

(iii) Define $\bW^\Gamma {\stackrel{\mathrm{def}}{=}}\{w\in \bW|\ell(ws)>\ell(w) \ \textrm{for  all} \ s\in\Gamma\}$. Given $w\in \bW$, there is a unique $u\in \bW^\Gamma$ and a unique $v\in \bW_\Gamma$ such that $w=uv$. Their lengths satisfy $\ell (w)=\ell(u)+\ell (v)$. Moreover, $u$ is the unique element of shortest length in the coset $w\bW_\Gamma$.
\end{prop}

For all $w \in \bW$ we set $w=w^{\Gamma}w_{\Gamma}$ with 
$w^{\Gamma} \in 
\bW^{\Gamma}$ and $w_{\Gamma} \in \bW_{\Gamma}$. Then if 
$\beta \in \bW_{\Gamma}$ one has 
$l(w\beta)=l(w^{\Gamma})+l(w_{\Gamma}\beta)$.

From (\ref{rel:bordoR}) it follows:
\begin{equation}\label{eq:bb}
\partial(E(w,\Gamma)) = 
w^{\Gamma}.\partial(E(w_{\Gamma},\Gamma)) 
\end{equation}
where the action (\ref{prodperaz}) is extended to $C(\bW)$ by linearity.

As a consequence we have a direct sum decomposition into isomorphic factors: 
\begin{equation}\label{similcomol}
H_q(G_n^k,R_{\tau}) \simeq \bigoplus_{j=1}^{\mid \bW^{S_k} \mid}H_q(C(\bW_{S_k}),R_{\tau}).
\end{equation}

\section{Preparation for the Main Theorem} 

Let $m_k:=\mid \bW^{S_k} \mid$ and $\bW_k:=\bW_{S_k}$; the exact sequences 
(\ref{eqn:suc.es.}) with relations (\ref{similcomol}) 
give rise to the corresponding long 
exact sequences in homology
\begin{equation}\label{cohom} 
\begin{split}
\cdots \lgr H_{q+1}(F^{k+1}_n(\bW),R_{\tau}) 
&\lgr \bigoplus_{j=1}^{m_{n-k-1}}H_{q-k}(C(\bW_{S_{n-k-1}}),R_{\tau}) 
\lgr\\
&\lgr H_{q}(F^{k}_n(\bW),R_{\tau}) \lgr H_{q}(F^{k+1}_n(\bW),R_{\tau}) \lgr 
\cdots .
\end{split} 
\end{equation}
 
We have the following

\begin{lemma} \label{stabile}
If $H_{q-h}(C(\bW_{n-h-1}),R_{\tau})$ are trivial $\Z$-modules
for all h such that $k \leq h \leq q$, then 
$H_{q}(F^{k}_n(\bW),R_{\tau})$ is also trivial.
\end{lemma}

\textbf{Proof:} 
From (\ref{eqn:suc.es.}) and (\ref{similcomol}) one has the exact sequences of complexes
\begin{equation} \label{eqn2}
\begin{split}
&0 \lgr \bigoplus_{j=1}^{m_{n-k-1}}C(\bW_{n-k-1})[k] \lgr 
F^k_n(\bW) \lgr F^{k+1}_{n}(\bW) \lgr 0\\
&0 \lgr \bigoplus_{j=1}^{m_{n-k-2}}C(\bW_{n-k-2})[k+1] \lgr 
F^{k+1}_n(\bW) \lgr F^{k+2}_{n}(\bW) \lgr 0\\
&\cdots \\
&0 \lgr \bigoplus_{j=1}^{m_{n-q-1}}C(\bW_{n-q-1})[q] \lgr 
F^{q}_n(\bW) \lgr F^{q+1}_{n}(\bW) \lgr 0\\
\end{split} 
\end{equation}   

The last sequence gives rise to the long exact sequence in homology:
\begin{equation} \label{eqnn1}
\cdots \lgr \bigoplus_{j=1}^{m_{n-q-1}}H_0(C(\bW_{n-q-1}),R_{\tau})
\lgr H_q(F^q_{n}(\bW),R_{\tau}) \lgr 0.
\end{equation}

By hypothesis $H_0(C(\bW_{n-q-1}),R_{\tau})$ is a trivial $\Z$-module then
$H_q(F^q_{n},R_{\tau})$ is also trivial.
 
We get the thesis going backwards in (\ref{eqn2}) and considering, in a 
similar way of (\ref{eqnn1}), the long 
exact sequences induced. $\qed$\\

Recall (see (\ref{spezz})) the decomposition:
$$
H_*(C(\bW),R_{\tau}) = \bigoplus_{r|\sharp(\A(\bW))} [R/ (\varphi_r)]^{\alpha_r}.
$$

It follows that if 
$\sharp(\A(\bW))$ and $\sharp(\A(\bW_{n-h}))$ are coprimes, the maps $i_{n,h}$ 
of (\ref{eqn3}) give rise to homology maps with images
sums of copies of $\{\varphi_1\}$
($\{\varphi_1\}^{0}$ means that the map 
is identically $0$).

We have that $\sharp(\A(\aemme))=n(n+1)/2$ and 
$\sharp(\A(\bemme))=n^2$ (see \cite{Bou}). If we fix
\begin{equation*}
\begin{split} 
&(n,h)=(3q+1,2) \mbox{ for } \aemme \\
&(n,h)=(n,1)    \mbox{ for } \bemme
\end{split}
\end{equation*}
then 
\begin{equation*}
\begin{split} 
&(\sharp(\A(\mathbf{A_{3q+1}})),\sharp(\A(\mathbf{A_{3q-1}})))=1\\
&(\sharp(\A(\bemme)),\sharp(\A(\mathbf{B_{n-1}})))=1 .
\end{split}
\end{equation*}

Since $i_{n,h}$ are injective, we can complete
(\ref{eqn3}) to short exact sequences of complexes which give, by the 
above remark:
\begin{equation} \label{eqn4}
\begin{split}
0 \lgr \oplus \{\varphi_1\} \lgr H_q(C(\mathbf{A_{3q+1}}),&R_{\tau}) \lgr H_q(C(\mathbf{A_{3q+1}})/\oplus_{j=1}^{m_{3q-1}} C(\mathbf{A_{3q-1}}),R_{\tau}) \lgr \\ 
&\bigoplus_{j=1}^{m_{3q-1}} H_{q-1}(C(\mathbf{A_{3q-1}}),R_{\tau}) \lgr 
\oplus \{\varphi_1\} \lgr \cdots
\end{split} 
\end{equation}
in case $\aemme$ and
\begin{equation} \label{eqn5}
\begin{split}
0 \lgr \oplus \{\varphi_1\} \lgr H_q(C(\bemme),&R_{\tau}) \lgr 
H_q(C(\bemme)/\oplus_{j=1}^{m_{n-1}} C(\mathbf{B_{n-1}}),R_{\tau}) \lgr \\ 
&\bigoplus_{j=1}^{m_{n-1}} H_{q-1}(C(\mathbf{B_{n-1}}),R_{\tau}) \lgr 
\oplus \{\varphi_1\} \lgr \cdots .
\end{split} 
\end{equation}
in case $\bemme$.

In order to prove theorem \ref{teor:stab}, we need to study the complexes
$C(\mathbf{A_{3q+1}})/\oplus_{j=1}^{m_{3q-1}}~C(\mathbf{A_{3q-1}})$ and 
$C(\bemme)/\oplus_{j=1}^{m_{n-1}} C(\mathbf{B_{n-1}})$.

The latter is exactly the complex $F^1_{n}(\bemme)$.

The farmer is the complex with basis over $R$: 
$$
\mathcal{E}_T:=\{E(w,\Gamma \cup T) \mid
w \in \mathbf{A_{3q+1}} \mbox{ and } \Gamma \subset S_{3q-1}\}
$$
for $\emptyset \subsetneq T \subset S^{3q-1}$.
We remark that $\mathcal{E}_{\{s_{3q}\}}$
is the basis of a complex
isomorphic to $(3q+2)$ copies of $F^1_{3q}(\mathbf{A_{3q}})$,  
$\mathcal{E}_{\{s_{3q+1}\}}$ generates
the subcomplex given by the image
of $G_{3q+1}^{3q-1}(\mathbf{A_{3q+1}})$ 
by the map $i_{3q+1}[1]$ and the elements of 
$\mathcal{E}_{\{s_{3q+1},s_{3q}\}}$ 
are the generators of the module $F^2_{3q+1}(\mathbf{A_{3q+1}})$.

Now we set
$$
(F^k_n(\bW))_h:=\{E(w,\Gamma) \in F^k_n(\bW) \mid \mid \Gamma \mid =h\}
$$
and $\partial^k_{n,h}: (F^k_n(\bW))_h \lgr (F^k_n(\bW))_{h-1}$ the $h$-th 
boundary
map in $F^k_n(\bW)$ ($\partial_{n,h}:=\partial^0_{n,h}$ is the boundary
map in $C(\bW)_h$).

Then the $h$-th boundary matrix of 
$C(\mathbf{A_{3q+1}})/\oplus_{j=1}^{m_{3q-1}} C(\mathbf{A_{3q-1}})$ is of 
the form
$$
{\overline{\partial}_h} = \left[
\begin{array}{ccc}
\oplus_{i=1}^{3q+2}\partial^1_{3q,h} & 0                 & A_1 \\
    0 & \oplus_{i=1}^{\frac{(3q+1)(3q+2)}{2}}\partial_{3q-1,h-1} & A_2 \\
          0             & 0                             & \partial^2_{3q+1,h} 
\end{array}
\right]
$$ 
where
$A_1$ and $A_2$ are the matrices of the image of the generators in  
$\mathcal{E}_{\{s_{3q},s_{3q+1}\}}$ restricted to   
$\mathcal{E}_{\{s_{3q}\}}$ and $\mathcal{E}_{\{s_{3q+1}\}}$ respectively.

Moreover all homology groups of the complexes 
$F^k_n(\bW)$ are torsion groups so the rank of $\partial^k_{n,h}$ equals the
rank of ker$(\partial^k_{n,h-1})$. Then it is not difficult to see that
the rank of $\overline{\partial}_{h}$ is exactly the sum of  
$(3q+2)$ times the rank of $\partial^1_{3q,h}$, 
$\frac{(3q+1)(3q+2)}{2}$ times the rank of $\partial_{3q-1,h-1}$ and 
the rank of $\partial^2_{3q+1,h}$.

\begin{rem}\label{oss:quoz}It follows that in order to prove that 
$H_k(C(\mathbf{A_{3q+1}})/~\oplus_{j=1}^{m_{3q-1}}~C(\mathbf{A_{3q-1}}),R_{\tau})$
is sum of copies of $\{\varphi_1\}$, i.e. a trivial $\Z$-module,
it is sufficient to prove the same result for 
$H_k(F^1_{3q}(\mathbf{A_{3q}}),R_{\tau})$, $H_{k-1}(C(\mathbf{A_{3q-1}}),R_{\tau})$ and 
$H_k(F^2_{3q+1}(\mathbf{A_{3q+1}}),R_{\tau})$.
\end{rem}

\section{Proof of the Main Theorem} 

In this section we prove theorem \ref{teor:stab}. 
This is equivalent to prove that 
$H_k(C(\aemme),\R_{\tau})$ is a trivial $\Z$-module for $n \geq 3k+1$, 
$H_k(C(\bemme),\R_{\tau})$ is trivial for $n \geq 2k+1$ and
$H_k(C(\demme),\R_{\tau})$ is trivial for $n \geq 3k+2$ (see relation (\ref{comom})).

For cases $\aemme$ and $\bemme$
 we use induction on the degree of homology. Case $\demme$ will follow
from $\aemme$.

By standard methods (see also \cite{Simona}) one gets the first step of induction, which  is 
\begin{equation}\label{passo0}
H_0(C(\aemme),R_{\tau}) \simeq H_0(C(\bemme),R_{\tau}) \simeq \{\varphi_1\} 
\end{equation}
for all $n \geq 1$.

One supposes that
$H_{k-1}(C(\aemme),R_{\tau})$ and 
$H_{k-1}(C(\bemme),R_{\tau})$ are trivial $\Z$-modules, respectively,
for all $n \geq 3(k-1)+1$ and $n \geq 2(k-1)+1$.

We have to prove that $H_k(C(\aemme),R_{\tau})$ and 
$H_k(C(\bemme),R_{\tau})$ are trivial $\Z$-modules, respectively,
for all $n \geq 3k+1$ and $n \geq 2k+1$.

First we consider the case $n=3k+1$ ( $n=2k+1$ ); using the sequence 
(\ref{eqn4}) ( (\ref{eqn5}) ), one needs only to prove 
that $H_k(C(\mathbf{A_{3k+1}})/\oplus_{j=1}^{m_{3k-1}} 
C(\mathbf{A_{3k-1}}),R_{\tau})$\\
( $H_k(C(\mathbf{B_{2k+1}})/\oplus_{j=1}^{m_{2k}}~C(\mathbf{B_{2k}}),R_{\tau})$ ) is trivial.

The assertion in case $\mathbf{B_{2k+1}}$ follows from Lemma \ref{stabile} 
since
$$
H_*(C(\mathbf{B_{2k+1}})/\oplus_{j=1}^{m_{2k}} C(\mathbf{B_{2k}}),R_{\tau})= 
H_*(F^1_{2k+1}(\mathbf{B_{2k+1}}),R_{\tau})
$$
and $H_{k-h}(C(\mathbf{B_{2k-h}}),R_{\tau})$ is trivial for all 
$ 1 \leq h \leq k$ by inductive hypothesis.

The proof in case $\mathbf{A_{3k+1}}$ is a consequence of 
remark \ref{oss:quoz}. 

One has that
$H_{k-1}(C(\mathbf{A_{3k-1}}),R_{\tau})$ is a trivial $\Z$-module
by induction and, from Lemma \ref{stabile}, 
$H_k(F^1_{3k}(\mathbf{A_{3k}}),R_{\tau})$ and  
$H_k(F^2_{3k+1}(\mathbf{A_{3k+1}}),R_{\tau})$ are trivial since  
$H_{k-h}(C(\mathbf{A_{3k-h-1}}),R_{\tau})$
and 
$H_{k-h}(C(\mathbf{A_{3k-h}}),R_{\tau})$ are trivial 
by hypothesis, respectively, for $1 \leq h \leq k$ and $2 \leq h \leq k$.\\

Let now $n > 3k+1$, we conclude the proof for $\aemme$
using induction on $n$. 
%
%

One supposes that $H_k(C(\mathbf{A_{n-1}}),R_{\tau})$ is trivial as 
$\Z$-module. Moreover
$H_{k-h}(C(\mathbf{A_{n-h-1}}),R_{\tau})$ are trivial by inductive 
hypothesis on the degree of homology, since
$(n-h-1)\geq~3(k-h)+1$ for all $1 \leq h \leq k$. 
Then $H_{k-h}(C(\mathbf{A_{n-h-1}}),R_{\tau})$ are trivial for 
$0 \leq h \leq k$ and the thesis follows from Lemma \ref{stabile}.

The proof in case $\bemme$, for $n> 2k+1$, is exactly the same.\\

Case $\demme$ is a consequence of Lemma \ref{stabile} applied to the 
exact sequence of complexes
$$
0 \lgr \bigoplus_{j=1}^{m_{n-1}}C(\mathbf{D_{S_{n-1}}}) \lgr C(\demme) 
\lgr F^{1}_n(\demme) \lgr 0
$$
since $C(\mathbf{D_{S_k}})=C(\mathbf{A_k})$ for all $0 \leq k \leq n-1$ 
(we use the standard Dynking diagram of $\demme$).
$\qed$\\

The last step is the\\
\\
\textbf{Proof of theorem \ref{teor:esp}}
From the universal coefficients theorem it follows
\begin{equation} \label{eqn8}
H_k(C(\bW),\{\varphi_1\}) \simeq H_k(C(\bW),R_{\tau}) \otimes 
\{\varphi_1\} \oplus Tor(H_{k-1}(C(\bW),R_{\tau}),\{\varphi_1\}).
\end{equation}

If we set
\begin{equation*}
rk_{\Q}(H_k(C(\bW),R_{\tau}) \otimes \{\varphi_1\})=:a_{k+1}
\end{equation*}
then, in the range specified in theorem \ref{teor:stab}
\begin{equation*}
rk_{\Q}[Tor(H_{k-1}(C(\bW),R_{\tau}),\{\varphi_1\})]=:a_{k}.
\end{equation*}

We recall, also, that $\{\varphi_1\}=\Q$, then
$$
H_k(C(\bW),\{\varphi_1\})=H_k(C(\bW),\Q),
$$
moreover the rank of $H_k(C(\bW),\Q)$ equals the rank of
$H^k(C(\bW),\Z)$.

It follows that relation (\ref{eqn8}) gives
$$
rk[H^k(C(\bW),\Z)]=a_{k+1}+a_{k}
$$
and from a simple induction
$$
a_{k+1}=\sum_{i=0}^{k}(-1)^{(k-i)} rk H^i(C(\bW),\Z). \qed
$$
                                    
\begin{rem} With the same technique used to prove theorem \ref{teor:stab},
it is possible to prove a more general result.

Let $(\bW,S)$ be a finite Coxeter system with $\mid S\mid = n$ and $m \in \N$ 
s.t. $m~\mid~o(\A(\bW))$. If there exists an integer $h$ s.t.
$m \nmid o(\A(\bW_k))$ for all $h < k < n$, 
then there exists an integer $p$ s.t., for all $r < p$,  
$H^r(C(\bW_h),R_{\tau})$ is annihilated by a squarefree element 
$(1-\tau^s)$ with $s\mid o(\A(\bW))$, $s < m$, and, for all $q < p+(n-~h~-~1)$,
$H^q(C(\bW),R_{\tau})$ 
is annihilated by a squarefree element 
$(1-\tau^a)$ with $a\mid o(\A(\bW))$, $a < m$.

As corollaries we obtain: 

\begin{itemize}
\item 
$H^{q+1}(C(\mathbf{A_{3q}}),R_{\tau})$ and
$H^{q+1}(C(\mathbf{A_{3q-1}}),R_{\tau})$ are annihilated by the squarefree
element $(1-\tau^3)$. 
\item if $m \mid o(\A(\bW))$ and 
$m \nmid o(\A(\bW_k))$ for all $k < n$ then, for $h < n$,
$H^h(C(\bW),R_{\tau})$ is annihilated by a squarefree element 
$(1-\tau^s)$ with $s\mid o(\A(\bW))$, $s < m$.
\end{itemize}
\end{rem}
\vspace{1,5cm}
\emph{\textbf{Acknowledgements} Special thanks to prof.Mario Salvetti for very 
useful conversations.}

\phantom{\cite{Bou,boss1,Simona}}
\phantom{\cite{boss4,boss3,Decsal,Denham}}
\phantom{\cite{b1,Brie,boss2,b5,OS}}
\phantom{\cite{b8,b9,b10,b11,b13,b14,b15}}
\newpage{\pagestyle{empty}\cleardoublepage}
\bibliographystyle{plain}
\addcontentsline{toc}{chapter}{Bibliografia}
\bibliography{maint} 

\end{document}